\documentclass[a4paper]{amsart}
\usepackage{graphicx,amsmath,amsfonts,latexsym,amssymb,amsthm,mathrsfs, color,hyperref}
\usepackage[latin1]{inputenc}
\usepackage{amsthm}
\usepackage{amssymb}
\theoremstyle{plain}
\newtheorem{thm}{Theorem}[section]

\newtheorem{corollary}[thm]{Corollary}
\newtheorem{definition}[thm]{Definition}
\newtheorem{theorem}[thm]{Theorem}

\newtheorem*{prooft1*}{Proof of Theorem 1.1}
\newtheorem*{proofc1*}{Proof of Corollary 1.2}

\begin{document}
\title[Conic manifolds under the Yamabe flow]%
{Conic manifolds under the Yamabe flow}
\author{Nikolaos Roidos}
\address{Department of Mathematics, University of Patras, 26504 Rio Patras, Greece}
\email{roidos@math.upatras.gr}

\thanks{The author was supported by Deutsche Forschungsgemeinschaft, grant SCHR 319/9-1.}
\date{\today}
\begin{abstract} 
We consider the unnormalized Yamabe flow on manifolds with conical singularities. Under certain geometric assumption on the initial cross-section we show well posedness of the short time solution in the $L^q$-setting. Moreover, we give a picture of the deformation of the conical tips under the flow by providing an asymptotic expansion of the evolving metric close to the boundary in terms of the initial local geometry. Due to the blow up of the scalar curvature close to the singularities we use maximal $L^q$-regularity theory for conically degenerate operators.
\end{abstract}
\maketitle

\section{Introduction}

Let $\mathcal{B}$ be a smooth compact $(n+1)$-dimensional manifold, $n\ge 2$, with closed (i.e. compact without boundary) possibly disconnected smooth boundary $\partial\mathcal{B}$ of dimension $n$. We endow $\mathcal{B}$ with a degenerate Riemannian metric $g_{0}=\{g_{ij}\}_{i,j\in\{1,...,n+1\}}$ which in a collar neighborhood $[0,1)\times\partial\mathcal{B}$ of the boundary is of the form
\begin{gather*}\label{metric}
g_{0}=dx^2+x^2h(x),
\end{gather*}
where $h(x)=\{h_{ij}(x)\}_{i,j\in\{1,...,n\}}$, $x\in[0,1)$, is a smooth up to $x=0$ family of non-degenerate up to $x=0$ Riemannian metrics on the cross section $\partial\mathcal{B}$. We call $\mathbb{B}=(\mathcal{B},g_{0})$ {\em manifold with conical singularities} or {\em conic manifold}; the singularities, i.e. the conical tips, correspond to the boundary $\{0\}\times \partial\mathcal{B}$ of $\mathcal{B}$.

We are interested in the unnormalized Yamabe flow on $\mathbb{B}$. Namely, we search for a family of Riemannian metrics $\{g(t)\}_{t\geq0}$ on $\mathcal{B}$ satisfying
\begin{eqnarray}\label{YFA}
g'(t)+R_{g(t)}g(t)&=&0, \quad t\in(0,T),\\\label{YFB}
g(0)&=&g_{0},
\end{eqnarray}
where $R_{g(t)}$ is the scalar curvature on $(\mathcal{B},g(t))$ and $T>0$.

We look as usual for solutions of \eqref{YFA}-\eqref{YFB} in the conformal class of $g_{0}$. The parabolic equation obtained, see \eqref{YF1}-\eqref{YF2}, has a quasilinear term that is described by the degenerate Laplacian on $\mathbb{B}$ and a {\em forcing term} involving the scalar curvature of $\mathbb{B}$, which blows up close to the conical tips. We regard the Laplacian as a {\em cone differential operator} acting on weighted {\em Mellin-Sobolev spaces} and employ the related maximal $L^{q}$-regularity theory, initially developed in \cite{CSS0}, \cite{CSS1}. Then, existence, uniqueness and maximal $L^q$-regularity of the short time solution is obtained by a theorem of Cl\'ement and Li. Moreover, we show that this solution becomes instantaneously smooth in space. In order to include the scalar curvature of $\mathbb{B}$ into our underline $L^p$-space we impose the following geometric assumption on the cross-section of $\mathbb{B}$, namely
\begin{gather}\label{assumptionR}
\bigg\{\begin{array}{cc}\text{Assume that the scalar curvature $R_{h(0)}$ of the cross-section $\partial\mathbb{B}=(\partial\mathcal{B},h(0))$ satisfies:}\\ R_{h(0)}=n(n-1). \end{array}
\end{gather}
Concerning the relation of the above condition to the solvability of the Yamabe problem in an asymptotic class of singular metrics we refer to \cite{AkBo}. Denote by $\mathbb{R}_{\omega}^{+}$ the set of smooth functions on $\mathbb{B}$ that are locally strictly positive constants close to the boundary (see Definition \ref{constfunt}). Although under assumption \eqref{assumptionR} the initial scalar curvature still blows up (see \eqref{scalcur}) we show the following.

\begin{theorem}\label{tshortyf}
Assume that \eqref{assumptionR} is satisfied and denote 
$$ 
\gamma_{0}=\min\Big\{-1+\sqrt{\Big(\frac{n-1}{2}\Big)^{2}-\lambda_{1}},\frac{n-1}{2}\Big\},
$$
where $\lambda_{1}$ is the greatest non-zero eigenvalue of the Laplacian on $\partial\mathcal{B}$ induced by the metric $h(0)$. Then, for any
\begin{gather}\label{datayf}
s\geq0, \quad \delta\in\Big(0,\gamma_{0}-\frac{n-3}{2}-\frac{2}{q}\Big) \quad \text{and} \quad p,q\in(1,\infty) \quad\text{with} \quad\frac{n+1}{p}+\frac{2}{q}<1,
\end{gather}
there exists a $T>0$ and a unique 
\begin{eqnarray*}\nonumber
\lefteqn{f\in W^{1,q}([0,T];\mathcal{H}_{p}^{s,\gamma_{0}-\delta}(\mathbb{B}))\cap L^{q}([0,T];\mathcal{H}_{p}^{s+2,\gamma_{0}+2-\delta}(\mathbb{B})\oplus\mathbb{R}_{\omega}^{+})}\\\label{solone}
&&\hookrightarrow \bigcap_{\varepsilon>0}C([0,T];\mathcal{H}_{p}^{s+2-\frac{2}{q}-\varepsilon,\gamma_{0}+2-\frac{2}{q}-\delta-\varepsilon}(\mathbb{B})\oplus\mathbb{R}_{\omega}^{+})\hookrightarrow C([0,T];C(\mathbb{B}))
\end{eqnarray*}
such that $g(t)=f(t)g_{0}$ solves \eqref{YFA}-\eqref{YFB}. Moreover, for any $\tau\in(0,T)$ we have that
\begin{eqnarray*}\nonumber
\lefteqn{f\in \bigcap_{\nu\geq0}W^{1,q}([\tau,T];\mathcal{H}_{p}^{\nu,\gamma_{0}-\delta}(\mathbb{B}))\cap L^{q}([\tau,T];\mathcal{H}_{p}^{\nu,\gamma_{0}+2-\delta}(\mathbb{B})\oplus\mathbb{R}_{\omega}^{+})}\\\label{soltwo}
&&\hookrightarrow \bigcap_{\nu\geq0,\varepsilon>0}C([\tau,T];\mathcal{H}_{p}^{\nu,\gamma_{0}+2-\frac{2}{q}-\delta-\varepsilon}(\mathbb{B})\oplus\mathbb{R}_{\omega}^{+}).
\end{eqnarray*}
\end{theorem}

The above result together with standard properties of Mellin-Sobolev spaces provide a picture of the deformation of the conical tips under the flow as well as its relation to the initial local geometry.

\begin{corollary}
Let $s$, $\delta$, $p$, $q$, $\lambda_{1}$ be as in \eqref{datayf}, $\varepsilon\in\big(0,2-\frac{2}{q}-\frac{n+1}{p}\big)$ and $T$, $f$ be as in Theorem \ref{tshortyf}. Then, there exist 
$$
f_{1}\in C([0,T];\mathbb{R}_{\omega}^{+}) \quad \text{and} \quad f_{2}\in C([0,T];\mathcal{H}_{p}^{s+2-\frac{2}{q}-\varepsilon,\gamma_{0}+2-\frac{2}{q}-\delta-\varepsilon}(\mathbb{B}))\hookrightarrow C([0,T];C(\mathbb{B}))
$$ 
with $f_{1}(0)=1$ and $f_{2}(0)=0$ such that in local coordinates $(x,y)$ on the collar part $(0,1)\times\partial\mathcal{B}$ we have
$$
g(t)=(f_{1}(t)+f_{2}(t,x,y))(dx^{2}+x^{2}h(x)),\quad t\in[0,T], \quad \text{with}\quad |f_{2}(t,x,y)|\leq Cx^{\beta-\varepsilon},
$$
where 
$$
\beta=\min\Big\{-\frac{n-1}{2}+\sqrt{\Big(\frac{n-1}{2}\Big)^{2}-\lambda_{1}},1\Big\}-\frac{2}{q}-\delta>0
$$ 
and the constant $C>0$ depends only on the data $\{\mathbb{B},\delta,p,q,\varepsilon,T\}$. 
\end{corollary}

In the case of a closed manifold the Yamabe flow as well as its normalized version (see, e.g., \cite[(5)]{Bre}) have been extensively studied in many aspects (existence, regularity, convergence, etc.). Concerning this situation, and in order to avoid the large amount of literature, for a complete introduction to the problem we only refer to the survey \cite{Bre} and to the references therein. 

The main difficulty of \eqref{YFA}-\eqref{YFB} compared to the classical case is the degeneracy of the associated Laplacian and at the same time the blow up of the scalar curvature close to the conical tips. In general, the Yamabe flow on manifolds with singularities is less studied. We close this section by pointing out some important contributions to this direction. The results stated below agree with ours in the case of conic manifolds. However, in our situation the geometric assumption is more flexible (see, e.g., \cite[(1.2)]{BaVer2}), $L^{q}$-regularity is established and moreover asymptotic behavior of the solution is provided and its interplay with the initial local geometry is shown.

In \cite{Sha1} and \cite{Sha2} the problem is considered on {\em singular manifolds} (including conic manifolds) in the sense of Amann (see, e.g., \cite{Am2}, \cite{Am1}, \cite{Am3}) by using maximal continuous and maximal $L^q$-regularity theory, respectively. By showing maximal regularity results for certain classes of degenerate differential operators (see \cite[Theorem 3.5]{Sha1}, \cite[Theorem 3.8]{Sha1}, \cite[Theorem 4.9]{Sha2} as well as \cite{Sha4} and \cite{Sha3}), in \cite[Theorem 4.4]{Sha1} and \cite[Theorem 5.1]{Sha2} it is shown under some assumptions that the type of the singularities is preserved under the evolution for short time. 

In \cite{BaVer} and \cite{BaVer2} the Yamabe flow is studied on manifolds with edge type singularities. By obtaining heat kernel estimates and then establishing mapping properties of the heat operator on appropriate H\"older spaces, under certain geometric assumptions in \cite[Theorem 1.7]{BaVer2} it is shown short time existence of a smooth solution which is improved in \cite[Theorem 1.3]{BaVer} to a long time existence and convergence result for the normalized Yamabe flow. For similar problems on conic manifolds we also refer to \cite{AkBo}, \cite{Ak}, \cite{JeRo}, \cite{KrVe}, \cite{Mazz} and \cite{Ver}.

\section{Maximal $L^q$-regularity for quasilinear parabolic equations}
\setcounter{equation}{0}

Let $X_{1}\overset{d}{\hookrightarrow} X_{0}$ be a continuously and densely injected complex Banach couple.
\begin{definition}[Sectorial operators]
Let $\mathcal{P}(K,\theta)$, $\theta\in[0,\pi)$, $K\geq1$, be the class of all closed densely defined linear operators $A$ in $X_{0}$ such that 
$$
S_{\theta}=\{\lambda\in\mathbb{C}\,|\, |\arg(\lambda)|\leq\theta\}\cup\{0\}\subset\rho{(-A)} \quad \mbox{and} \quad (1+|\lambda|)\|(A+\lambda)^{-1}\|_{\mathcal{L}(X_{0})}\leq K \quad \text{when} \quad \lambda\in S_{\theta}.
$$
The elements in $\mathcal{P}(\theta)=\cup_{K\geq1}\mathcal{P}(K,\theta)$ are called {\em (invertible) sectorial operators of angle $\theta$}.
\end{definition}

Consider the following abstract parabolic first order Cauchy problem
\begin{eqnarray}\label{app1}
u'(t)+Au(t)&=&w(t), \quad t\in(0,T),\\\label{app2}
u(0)&=&0,
\end{eqnarray}
where $-A:X_{1}\rightarrow X_{0}$ is the infinitesimal generator of an analytic semigroup on $X_{0}$ and $w\in L^q([0,T];X_{0})$, $q\in(1,\infty)$, $T>0$. The operator $A$ has {\em maximal $L^q$-regularity} if for any $w\in L^q([0,T];X_{0})$ there exists a unique $u\in W^{1,q}([0,T];X_{0})\cap L^{q}([0,T];X_{1})$ solving \eqref{app1}-\eqref{app2}. In this case $u$ depends continuously on $w$. Furthermore, the above property is independent of $q$ and $T$. 

We recall next the following boundedness condition for the resolvent of an operator.
\begin{definition}[$R$-sectorial operators]
Denote by $\mathcal{R}(K,\theta)$, $\theta\in[0,\pi)$, $K\geq1$, the class of all operators $A\in \mathcal{P}(\theta)$ in $X_{0}$ such that for any choice of $\lambda_{1},...,\lambda_{N}\in S_{\theta}\backslash\{0\}$ and $x_{1},...,x_{N}\in X_0$, $N\in\mathbb{N}\backslash\{0\}$, we have
\begin{eqnarray*}
\Big\|\sum_{k=1}^{N}\epsilon_{k}\lambda_{k}(A+\lambda_{k})^{-1}x_{k}\Big\|_{L^{2}([0,1];X_0)} \leq K \Big\|\sum_{k=1}^{N}\epsilon_{k}x_{k}\Big\|_{L^{2}([0,1];X_0)},
\end{eqnarray*}
where $\{\epsilon_{k}\}_{k=1}^{\infty}$ is the sequence of the Rademacher functions. The elements in $\mathcal{R}(\theta)=\cup_{K\geq1}\mathcal{R}(K,\theta)$ are called {\em $R$-sectorial operators of angle $\theta$}. 
\end{definition}

If we restrict to the class of UMD (unconditionality of martingale differences property, see, e.g., \cite[Section III.4.4]{Am}) Banach spaces, then the following result holds.

\begin{theorem}[{\rm Kalton and Weis, \cite[Theorem 6.5]{KW1}}]
In a UMD Banach space any $R$-sectorial operator of angle greater than $\frac{\pi}{2}$ has maximal $L^{q}$-regularity. 
\end{theorem}

Let $q\in(1,\infty)$, $U$ be an open subset of $(X_{1},X_{0})_{\frac{1}{q},q}$, $A(\cdot): U\rightarrow \mathcal{L}(X_{1},X_{0})$ and $F(\cdot,\cdot): U\times [0,T_{0}]\rightarrow X_{0}$, for some $T_{0}>0$. Consider the problem
\begin{eqnarray}\label{aqpp1}
u'(t)+A(u(t))u(t)&=&F(u(t),t)+G(t),\quad t\in(0,T),\\\label{aqpp2}
u(0)&=&u_{0},
\end{eqnarray}
where $T\in(0,T_{0})$, $u_{0}\in U$ and $G\in L^{q}([0,T_{0}];X_{0})$. Maximal $L^q$-regularity for the linearization $A(u_{0})$ together with appropriate Lipschitz continuity conditions imply existence of a short time solution for the above equation.
\begin{theorem}[{\rm Cl\'ement and Li, \cite[Theorem 2.1]{CL}}]\label{ClementLi}
Assume that:
\begin{itemize}
\item[(H1)] $A(\cdot)\in C^{1-}(U;\mathcal{L}(X_{1},X_{0}))$.
\item[(H2)] $F(\cdot,\cdot)\in C^{1-,1-}(U\times [0,T_{0}];X_{0})$.
\item[(H3)] $A(u_{0})$ has maximal $L^q$-regularity.
\end{itemize}
Then, there exists a $T\in(0,T_{0})$ and a unique 
$$
u\in W^{1,q}([0,T];X_{0})\cap L^{q}([0,T];X_{1})
$$
solving \eqref{aqpp1}-\eqref{aqpp2}. 
\end{theorem}

Finally, recall the following embedding of the maximal $L^q$-regularity space, namely
\begin{gather}\label{interpemb}
W^{1,q}([0,T];X_{0})\cap L^{q}([0,T];X_{1})\hookrightarrow C([0,T];(X_{1},X_{0})_{\frac{1}{q},q}), \quad T>0, \,\, q\in(1,\infty),
\end{gather}
see e.g. \cite[Theorem III.4.10.2]{Am}.

\section{The Laplacian as a cone differential operator}

When the Laplacian $\Delta$ on $\mathbb{B}$ is restricted on the collar part $(0,1)\times\partial\mathcal{B}$ it takes the degenerate form
\begin{gather*}\label{delta}
\Delta=\frac{1}{x^{2}}\Big[(x\partial_{x})^{2}+\Big(n-1+\frac{x\partial_{x}\det[h(x)]}{2\det[h(x)]}\Big)(x\partial_{x})+\Delta_{h(x)}\Big],
\end{gather*}
where $\Delta_{h(x)}$ is the Laplacian on $\partial\mathcal{B}$ induced by the metric $h(x)$. We regard $\Delta$ as an element in the class of {\em cone differential operators} or {\em Fuchs type operators}, i.e. the naturally appearing degenerate differential operators on $\mathbb{B}$. In this section we recall some basic facts and results from the related underline pseudodifferential theory, i.e. the {\em cone calculus}, towards the direction of the study of nonlinear partial differential equations. For more details we refer to \cite{CSS0}, \cite{CSS1}, \cite{Le}, \cite{GKM}, \cite{RS2}, \cite{RS3}, \cite{RS1}, \cite{SS}, \cite{Schu} and \cite{Sei}.

We call cone differential operator of order $\mu\in\mathbb{N}$ any $\mu$-th order differential operator $A$ with smooth coefficients in the interior $\mathbb{B}^{\circ}$ of $\mathbb{B}$ such that when it is restricted to the collar part $(0,1)\times\partial\mathcal{B}$ it admits the following form 
\begin{gather}\label{Acone}
A=x^{-\mu}\sum_{k=0}^{\mu}a_{k}(x)(-x\partial_{x})^{k}, \quad \mbox{where} \quad a_{k}\in C^{\infty}([0,1);\mathrm{Diff}^{\mu-k}(\partial\mathbb{B})).
\end{gather}
If, in addition to the usual pseudodifferential symbol, we assume that the {\em rescaled symbol} (see e.g. \cite[(2.3)]{CSS1} for definition) of $A$ is also pointwise invertible, then $A$ is called {\em $\mathbb{B}$-elliptic}; this is the case for the Laplacian $\Delta$.

Cone differential operators act naturally on scales of {\em Mellin-Sobolev} spaces. Let $\omega\in C^{\infty}(\mathbb{B})$ be a fixed cut-off function near the boundary, i.e. a smooth non-negative function on $\mathcal{B}$ with $\omega=1$ near $\{0\}\times\partial \mathcal{B}$ and $\omega=0$ on $\mathcal{B}\backslash([0,1)\times \partial \mathcal{B})$. Moreover, assume that in local coordinates $(x,y)\in [0,1)\times \partial\mathcal{B}$, $\omega$ depends only on $x$. Denote by $C_{c}^{\infty}$ the space of smooth compactly supported functions and by $H_{p}^{s}$, $p\in(1,\infty)$, $s\in\mathbb{R}$, the usual Sobolev space modeled on $L^{p}$. 

\begin{definition}[Mellin-Sobolev spaces]
For any $\gamma\in\mathbb{R}$ consider the map 
$$
M_{\gamma}: C_{c}^{\infty}(\mathbb{R}_{+}\times\mathbb{R}^{n})\rightarrow C_{c}^{\infty}(\mathbb{R}^{n+1}) \quad \mbox{defined by} \quad u(x,y)\mapsto e^{(\gamma-\frac{n+1}{2})x}u(e^{-x},y). 
$$
Furthermore, take a covering $\kappa_{i}:U_{i}\subseteq\partial\mathcal{B} \rightarrow\mathbb{R}^{n}$, $i\in\{1,...,N\}$, $N\in\mathbb{N}\backslash\{0\}$, of $\partial\mathcal{B}$ by coordinate charts and let $\{\phi_{i}\}_{i\in\{1,...,N\}}$ be a subordinated partition of unity. For any $s\in\mathbb{R}$ and $p\in(1,\infty)$ let $\mathcal{H}^{s,\gamma}_p(\mathbb{B})$ be the space of all distributions $u$ on $\mathbb{B}^{\circ}$ such that 
$$
\|u\|_{\mathcal{H}^{s,\gamma}_p(\mathbb{B})}=\sum_{i=1}^{N}\|M_{\gamma}(1\otimes \kappa_{i})_{\ast}(\omega\phi_{i} u)\|_{H^{s}_p(\mathbb{R}^{n+1})}+\|(1-\omega)u\|_{H^{s}_p(\mathbb{B})}
$$
is defined and finite, where $\ast$ refers to the push-forward of distributions. The space $\mathcal{H}^{s,\gamma}_p(\mathbb{B})$, called {\em (weighted) Mellin-Sobolev space}, is independent of the choice of the cut-off function $\omega$, the covering $\{\kappa_{i}\}_{i\in\{1,...,N\}}$ and the partition $\{\phi_{i}\}_{i\in\{1,...,N\}}$; if $A$ is as in \eqref{Acone}, then it induces a bounded map
$$
A: \mathcal{H}^{s+\mu,\gamma+\mu}_p(\mathbb{B}) \rightarrow \mathcal{H}^{s,\gamma}_p(\mathbb{B}).
$$
Finally, if $s\in \mathbb{N}$, then equivalently, $\mathcal{H}^{s,\gamma}_p(\mathbb{B})$ is the space of all functions $u$ in $H^s_{p,loc}(\mathbb{B}^\circ)$ such that near the boundary
$$
x^{\frac{n+1}2-\gamma}(x\partial_x)^{k}\partial_y^{\alpha}(\omega(x) u(x,y)) \in L_{loc}^{p}\big([0,1)\times \partial \mathcal{B}, \sqrt{\mathrm{det}[h(x)]}\frac{dx}xdy\big),\quad k+|\alpha|\le s.
$$
\end{definition}

Let us restrict to the case of the Lapacian $\Delta$ and regard it as an unbounded operator in $\mathcal{H}^{s,\gamma}_p(\mathbb{B})$, $s,\gamma\in\mathbb{R}$, $p\in(1,\infty)$, with domain $C_{c}^{\infty}(\mathbb{B}^{\circ})$. The domain of its minimal extension (i.e. its closure) $\underline{\Delta}_{s,\min}$ is given by 
$$
\mathcal{D}(\underline{\Delta}_{s,\min})=\Big\{u\in \bigcap_{\varepsilon>0}\mathcal{H}^{s+2,\gamma+2-\varepsilon}_p(\mathbb{B}) \, |\, \Delta u\in \mathcal{H}^{s,\gamma}_p(\mathbb{B})\Big\}.
$$ 
If in addition the {\em conormal symbol} of $\Delta$, i.e. the following family of differential operators
$$
\lambda^{2}-(n-1)\lambda + \Delta_{h(0)} : \mathbb{C} \rightarrow \mathcal{L}(H_{2}^{2}(\partial\mathbb{B}),H_{2}^{0}(\partial\mathbb{B})),
$$
 is invertible on the line $\{\lambda\in\mathbb{C}\,|\, \mathrm{Re}(\lambda)= \frac{n-3}{2}-\gamma\}$, then we have precisely
$$
\mathcal{D}(\underline{\Delta}_{s,\min})=\mathcal{H}^{s+2,\gamma+2}_p(\mathbb{B}).
$$ 

The domain of the maximal extension $\underline{\Delta}_{s,\max}$ of $\Delta$, defined by
\begin{gather*}
\mathcal{D}(\underline{\Delta}_{s,\max})=\Big\{u\in\mathcal{H}^{s,\gamma}_p(\mathbb{B}) \, |\, \Delta u\in \mathcal{H}^{s,\gamma}_p(\mathbb{B})\Big\},
\end{gather*}
is expressed as
\begin{gather}\label{dmax1}
\mathcal{D}(\underline{\Delta}_{s,\max})=\mathcal{D}(\underline{\Delta}_{s,\min})\oplus\mathcal{E}_{\Delta,\gamma},
\end{gather}
where $\mathcal{E}_{\Delta,\gamma}$ is a finite-dimensional space called {\em asymptotics space}. $\mathcal{E}_{\Delta,\gamma}$ is independent of $s$; it consists of linear combinations of $C^{\infty}(\mathbb{B}^\circ)$ functions that vanish on $\mathcal{B}\backslash([0,1)\times\partial\mathcal{B})$ and in local coordinates $(x,y)\in (0,1)\times\partial\mathcal{B}$ they are of the form $\omega(x)c(y)x^{-\rho}\log^{m}(x)$ where $c\in C^{\infty}(\partial\mathbb{B})$, $\rho\in \mathbb{C}$ with 
$$
\frac{n-3}{2}-\gamma\leq\mathrm{Re}(\rho)<\frac{n-3}{2},
$$
and $m\in\{0,1\}$. The exponents $\rho$ are determined explicitly by the metric $h(\cdot)$. Due to \eqref{dmax1}, there are several closed extensions of $\Delta$ in $\mathcal{H}^{s,\gamma}_p(\mathbb{B})$ each one corresponding to a subspace of $\mathcal{E}_{\Delta,\gamma}$. For an overview on the domain structure of a general $\mathbb{B}$-elliptic cone differential operator we refer to \cite[Section 3]{GKM} or alternatively to \cite[Sections 2.2-2.3]{SS}.

\begin{definition}\label{constfunt}
Recall that $\partial\mathcal{B}=\cup_{i=1}^{k_{\mathcal{B}}}\partial\mathcal{B}_{i}$, for certain $k_{\mathcal{B}}\in\mathbb{N}\backslash\{0\}$, where $\partial\mathcal{B}_{i}$ are closed, smooth and connected. Denote by $\mathbb{C}_{\omega}$ the space of all $C^{\infty}(\mathbb{B}^\circ)$ functions $c$ that vanish on $\mathcal{B}\backslash([0,1)\times\partial\mathcal{B})$ and on each component $[0,1)\times\partial\mathcal{B}_{i}$, $i\in\{1,...,k_{\mathbb{B}}\}$, they are of the form $c_{i}\omega$, where $c_{i}\in\mathbb{C}$, i.e. $\mathbb{C}_{\omega}$ consists of smooth functions that are locally constant close to the boundary. Endow $\mathbb{C}_{\omega}$ with the norm $\|\cdot\|_{\mathbb{C}_{\omega}}$ given by $c\mapsto \|c\|_{\mathbb{C}_{\omega}}=(\sum_{i=1}^{k_{\mathcal{B}}}|c_{i}|^{2})^{\frac{1}{2}}$. Moreover, denote by $\mathbb{R}_{\omega}$ the subspace of $\mathbb{C}_{\omega}$ consisting of functions $c$ such that $c_{i}\in\mathbb{R}$ and by $\mathbb{R}_{\omega}^{+}$ its subset consisting of functions $c$ such that $c_{i}>0$.
\end{definition}

We close this section with a particular closed extension of the Laplacian. Under certain choice of the weight $\gamma$, $\mathbb{C}_{\omega}$ becomes a subspace of $\mathcal{E}_{\Delta,\gamma}$. By recalling that the Mellin-Sobolev spaces are UMD spaces, the corresponding realization of the Laplacian enjoys the property of maximal $L^q$-regularity as we can see from the following. 

\begin{theorem}[{\rm\cite[Theorem 3.2]{Ro2} or \cite[Theorem 5.6]{RS2}}]\label{RsecD}
Let $s\geq0$, $p\in(1,\infty)$ and 
\begin{eqnarray*}\label{choicegg}
\gamma\in\Big(\frac{n-3}2,\min\Big\{-1+\sqrt{\Big(\frac{n-1}{2}\Big)^{2}-\lambda_{1}} ,\frac{n+1}{2}\Big\}\Big),
\end{eqnarray*}
where $\lambda_{1}$ is the greatest non-zero eigenvalue of the boundary Laplacian $\Delta_{h(0)}$. Consider the closed extension $\underline{\Delta}_{s}$ of $\Delta$ in 
$$
X_{0}^{s}=\mathcal{H}_{p}^{s,\gamma}(\mathbb{B})
$$ 
with domain 
\begin{gather*}\label{DD}
\mathcal{D}(\underline{\Delta}_{s})=X_{1}^{s}=\mathcal{D}(\underline{\Delta}_{s,\min})\oplus\mathbb{C}_{\omega}=\mathcal{H}_{p}^{s+2,\gamma+2}(\mathbb{B})\oplus\mathbb{C}_{\omega}.
\end{gather*}
Then, for any $\theta\in[0,\pi)$ there exists some $c>0$ such that $c-\underline{\Delta}_{s}$ is $R$-sectorial of angle $\theta$.
\end{theorem}

\section{Conic manifolds under the Yamabe flow}

After setting $g(t)=u^{\frac{4}{n-1}}(t)g_{0}$, $t\geq0$, to \eqref{YFA}-\eqref{YFB} we obtain the following parabolic quasilinear equation over $u$, namely
\begin{eqnarray}\label{YF1}
u'(t)-nu^{-\frac{4}{n-1}}(t)\Delta u(t)&=&-\frac{n-1}{4}u^{\frac{n-5}{n-1}}(t)R_{g_{0}}, \quad t\in(0,T),\\\label{YF2}
u(0)&=&1,
\end{eqnarray}
where $R_{g_{0}}$ is the scalar curvature of $\mathbb{B}$. When $R_{g_{0}}$ is restricted on $(0,1)\times\partial\mathcal{B}$ it satisfies
\begin{gather}\label{scalcur}
xR_{g_{0}}=\frac{1}{x}(R_{h(0)}-n(n-1))+\text{non-singular terms},
\end{gather}
where $R_{h(0)}$ is the scalar curvature of $\partial\mathbb{B}$ (see also \cite[Theorem 2.1]{DoLa} or \cite[(2.8)]{JeRo}).

\subsection*{Proof of Theorem 1.1}

Let $\gamma=\gamma_{0}-\delta$ and $\nu\geq0$. By \cite[Lemma 3.2]{RS3} and \cite[Lemma 5.2]{RS3} we have the embedding
\begin{gather}\label{embd11}
\mathcal{H}_{p}^{\nu+2-\frac{2}{q}+\varepsilon,\gamma+2-\frac{2}{q}+\varepsilon}(\mathbb{B})\oplus\mathbb{C}_{\omega}\hookrightarrow (X_{1}^{\nu},X_{0}^{\nu})_{\frac{1}{q},q} \hookrightarrow \mathcal{H}_{p}^{\nu+2-\frac{2}{q}-\varepsilon,\gamma+2-\frac{2}{q}-\varepsilon}(\mathbb{B})\oplus\mathbb{C}_{\omega}\hookrightarrow C(\mathbb{B}),
\end{gather}
valid for all $\varepsilon>0$ sufficiently small. Let $U_{\nu}$ be an open bounded subset of $(X_{1}^{\nu},X_{0}^{\nu})_{\frac{1}{q},q}$ containing $1$. Due to \eqref{embd11}, we choose $U_{\nu}$ in such a way that there exists a finite closed path $\Gamma_{\nu}$ in $\{z\in\mathbb{C}\, |\, \mathrm{Re}(z)<0\}$ surrounding $\cup_{v\in U_{\nu}}\mathrm{Ran}(-v)$ with $\Gamma_{\nu}\cap\cup_{v\in U_{\nu}}\mathrm{Ran}(-v)=\emptyset$. For any $u_{1},u_{2}\in U_{\nu}$ and $\eta\in\mathbb{R}$ we have
\begin{gather}\label{eqinte3}
u_{1}^{\eta}-u_{2}^{\eta}=(u_{2}-u_{1})\frac{1}{2\pi i}\int_{\Gamma_{\nu}}(-z)^{\eta}(u_{1}+z)^{-1}(u_{2}+z)^{-1}dz.
\end{gather}
Moreover (cfr. \cite[Lemmas 3.2 and 6.2]{RS3}), the choice of the data implies the following:
\begin{gather}\label{BAlg}
\bigg\{\begin {array}{ll} \text{For all $\varepsilon>0$ sufficiently small, the space $\mathcal{H}_{p}^{\nu+2-\frac{2}{q}-\varepsilon,\gamma+2-\frac{2}{q}-\varepsilon}(\mathbb{B})\oplus\mathbb{C} _{\omega}$ is a Banach} \\
\text{algebra (up to norm equivalence) and also closed under holomorphic functional calculus.}
\end{array}
\end{gather}
Furthermore, by \eqref{embd11} and \cite[Lemma 6.3]{RS3}
\begin{gather}\label{eqinte345}
\bigg\{\begin {array}{ll} \text{For all $\varepsilon>0$ sufficiently small, the set $\{(v+z)^{-1}\,|\, v\in U_{\nu}, z\in\Gamma_{\nu}\}$}\\
\text{is bounded in $\mathcal{H}_{p}^{\nu+2-\frac{2}{q}-\varepsilon,\gamma+2-\frac{2}{q}-\varepsilon}(\mathbb{B})\oplus\mathbb{C}_{\omega}$.}
 \end{array}
\end{gather}
Finally, due to \cite[Corollary 3.3]{RS3} we have
\begin{gather}\label{BAlg22}
\bigg\{\begin{array}{ll} \text{For all $\varepsilon>0$ sufficiently small, elements in $\mathcal{H}_{p}^{\nu+2-\frac{2}{q}-\varepsilon,\gamma+2-\frac{2}{q}-\varepsilon}(\mathbb{B})\oplus\mathbb{C}_{\omega}$} \\
\text{act by multiplication as bounded operators on $X_{0}^{\nu}$.}
\end{array}
\end{gather}
We split the proof in several steps.

\subsubsection*{Short time existence} We apply Theorem \ref{ClementLi} to \eqref{YF1}-\eqref{YF2} with $X_{0}=X_{0}^{s}$, $X_{1}=X_{1}^{s}$,
$$
A(\cdot)=-n(\cdot)^{-\frac{4}{n-1}}\underline{\Delta}_{s}, \quad F(\cdot)=-\frac{n-1}{4}(\cdot)^{\frac{n-5}{n-1}}R_{g_{0}}, \quad G=0,
$$ 
$u_{0}=1$ and $U=U_{s}$. Concerning (H1), by \eqref{embd11}-\eqref{BAlg22} with $\nu=s$ we have that
\begin{eqnarray}\nonumber
\lefteqn{\|A(u_{1})-A(u_{2})\|_{\mathcal{L}(X_{1}^{s},X_{0}^{s})}}\\\nonumber
&=&n\|(u_{1}^{-\frac{4}{n-1}}-u_{2}^{-\frac{4}{n-1}})\underline{\Delta}_{s}\|_{\mathcal{L}(X_{1}^{s},X_{0}^{s})}\\\nonumber
&\leq&C_{1}\|u_{1}^{-\frac{4}{n-1}}-u_{2}^{-\frac{4}{n-1}}\|_{\mathcal{L}(X_{0}^{s})}\\\label{contnt}
&\leq&C_{2}\|u_{1}-u_{2}\|_{(X_{1}^{s},X_{0}^{s})_{\frac{1}{q},q}},
\end{eqnarray}
for certain constants $C_{1}, C_{2}>0$.

Similarly, concerning (H2), by noting that $R_{g_{0}}\in X_{0}^{s}$ due to \eqref{assumptionR}-\eqref{datayf}, for $\varepsilon>0$ sufficiently small we estimate
\begin{eqnarray*}
\lefteqn{\|F(u_{1})-F(u_{2})\|_{X_{0}^{s}}}\\
&\leq&C_{3}\|u_{1}^{-\frac{n-5}{n-1}}-u_{2}^{-\frac{n-5}{n-1}}\|_{\mathcal{H}_{p}^{s+2-\frac{2}{q}-\varepsilon,\gamma+2-\frac{2}{q}-\varepsilon}(\mathbb{B})\oplus\mathbb{C}_{\omega}}\|R_{g_{0}}\|_{X_{0}^{s}}\\
&\leq&C_{4}\|u_{1}-u_{2}\|_{(X_{1}^{s},X_{0}^{s})_{\frac{1}{q},q}},
\end{eqnarray*}
for some constants $C_{3}, C_{4}>0$.

The property (H3) follows immediately from Theorem \ref{RsecD}. Therefore, there exists a $T>0$ and a unique 
\begin{equation}\label{solone}
\begin{split}
u\in & \, W^{1,q}([0,T];X_{0}^{s}))\cap L^{q}([0,T];X_{1}^{s})\hookrightarrow C([0,T];(X_{1}^{s},X_{0}^{s})_{\frac{1}{q},q})\\
&\hookrightarrow \bigcap_{\varepsilon>0}C([0,T];\mathcal{H}_{p}^{s+2-\frac{2}{q}-\varepsilon,\gamma+2-\frac{2}{q}-\varepsilon}(\mathbb{B})\oplus\mathbb{C}_{\omega})\hookrightarrow C([0,T];C(\mathbb{B}))
\end{split}
\end{equation}
solving \eqref{YF1}-\eqref{YF2}, where for the last embedding we have used \eqref{interpemb} and \eqref{embd11}. By taking the complex conjugate in \eqref{YF1}-\eqref{YF2} and using the above uniqueness, we conclude that we can replace $\mathbb{C}_{\omega}$ in \eqref{solone} by $\mathbb{R}_{\omega}$.

\subsubsection*{Smoothness in space} We apply \cite[Theorem 3.1]{RS4} to \eqref{YF1}-\eqref{YF2} with $Y_{0}^{j}=X_{0}^{s+\frac{j}{q}}$, $Y_{1}^{j}=X_{1}^{s+\frac{j}{q}}$, $j\in\mathbb{N}$, $Z=\{v\in U_{s}\, |\, v\geq\alpha\}$, for certain $\alpha\in(0,1)$, and $A$, $F$ as above. According to \eqref{solone} we choose $T>0$ small enough such that $u(t)\in Z$, $t\in [0,T]$; in particular, we can replace $\mathbb{C}_{\omega}$ in \eqref{solone} by $\mathbb{R}_{\omega}^{+}$. 

Concerning the condition (i) of \cite[Theorem 3.1]{RS4}, by \eqref{BAlg}, \eqref{BAlg22} and \eqref{solone} we have that $A(u(t)):Y_{1}^{0}\rightarrow Y_{0}^{0}$, $t\in[0,T]$, is a well defined map, and moreover by \cite[Theorem 6.1]{RS3} it has maximal $L^{q}$-regularity. Furthermore, by \eqref{contnt}, we have
\begin{eqnarray*}
\|A(u(t_{1}))-A(u(t_{2}))\|_{\mathcal{L}(Y_{1}^{0},Y_{0}^{0})}\leq C_{5}\|u(t_{1})-u(t_{2})\|_{(Y_{1}^{0},Y_{0}^{0})_{\frac{1}{q},q}}, \quad t_{1},t_{2}\in[0,T],
\end{eqnarray*}
with some constant $C_{5}$ depending only on $u$ and $T$. Thus, $A(u(\cdot))\in C([0,T];\mathcal{L}(Y_{1}^{0},Y_{0}^{0}))$ due to \eqref{solone}.

Concerning the condition (ii) of \cite[Theorem 3.1]{RS4}, first note that \eqref{embd11} implies
\begin{gather}\label{embY}
(Y_{1}^{j},Y_{0}^{j})_{\frac{1}{q},q}\hookrightarrow\bigcap_{\varepsilon>0} \mathcal{H}_{p}^{s+2+\frac{j-2}{q}-\varepsilon,\gamma+2-\frac{2}{q}-\varepsilon}(\mathbb{B})\oplus\mathbb{C}_{\omega}, \quad j\in\mathbb{N}.
\end{gather}
Therefore, if $v\in Z\cap (Y_{1}^{j},Y_{0}^{j})_{\frac{1}{q},q}$, $j\in\mathbb{N}$, then by \eqref{BAlg}, \eqref{BAlg22} and \cite[Remark 2.8 (b)]{RS4} the operator $A(v):Y_{1}^{j+1}\rightarrow Y_{0}^{j+1}$ is well defined and moreover has maximal $L^{q}$-regularity. In addition, similarly to \eqref{contnt}, for any $w\in C([0,T];Z\cap (Y_{1}^{j},Y_{0}^{j})_{\frac{1}{q},q})$, $j\in\mathbb{N}$, and $t_{1},t_{2}\in[0,T]$ we have that
\begin{eqnarray*}\nonumber
\lefteqn{\|A(w(t_{1}))-A(w(t_{2}))\|_{\mathcal{L}(Y_{1}^{j+1},Y_{0}^{j+1})}}\\
&\leq&C_{6}\|w^{-\frac{4}{n-1}}(t_{1})-w^{-\frac{4}{n-1}}(t_{2})\|_{\mathcal{L}(Y_{0}^{j+1})}\\
&\leq&C_{7}\|w(t_{1})-w(t_{2})\|_{ \mathcal{H}_{p}^{s+2+\frac{j-2}{q}-\varepsilon,\gamma+2-\frac{2}{q}-\varepsilon}(\mathbb{B})\oplus\mathbb{C}_{\omega}},
\end{eqnarray*}
with some constants $C_{6},C_{7}>0$ depending only on $w$ and $T$. Due to \eqref{embY} we conclude $A(w(\cdot))\in C([0,T];\mathcal{L}(Y_{1}^{j+1},Y_{0}^{j+1}))$.

Finally, concerning the condition (iii) of \cite[Theorem 3.1]{RS4}, first note that $R_{g_{0}}\in \cap_{j\in\mathbb{N}}Y_{0}^{j}$. Moreover, by \eqref{eqinte3}-\eqref{eqinte345} and \eqref{embY} we have that
$$
w^{\frac{n-5}{n-1}}\in C([0,T]; \mathcal{H}_{p}^{s+2+\frac{j-2}{q}-\varepsilon,\gamma+2-\frac{2}{q}-\varepsilon}(\mathbb{B})\oplus\mathbb{C}_{\omega}) 
$$
for all $\varepsilon>0$. Therefore, by \cite[Corollary 3.3]{RS3} we deduce that $F(w(\cdot))\in L^{q}([0,T];Y_{0}^{j+1})$, $j\in\mathbb{N}$. We conclude that for each $\tau\in (0,T)$ we have
\begin{equation}\label{soltwo}
\begin{split}
u\in &\bigcap_{\nu\geq0}W^{1,q}([\tau,T];X_{0}^{\nu})\cap L^{q}([\tau,T];X_{1}^{\nu})\\
&\hookrightarrow \bigcap_{\nu\geq0}C([\tau,T];(X_{1}^{\nu},X_{0}^{\nu})_{\frac{1}{q},q}) \hookrightarrow \bigcap_{\nu\geq0,\varepsilon>0}C([\tau,T];\mathcal{H}_{p}^{\nu,\gamma+2-\frac{2}{q}-\varepsilon}(\mathbb{B})\oplus\mathbb{C}_{\omega}).
\end{split}
\end{equation}

\subsubsection*{Regularity of the conformal factor} Let $f=u^{\frac{4}{n-1}}$. By \eqref{eqinte3}-\eqref{BAlg22}, \eqref{solone} and \eqref{soltwo} for each $\eta\in\mathbb{R}$ we have 
\begin{equation}\label{ert}
\begin{split}
u^{\eta}\in & \bigcap_{\varepsilon>0}C([0,T];\mathcal{H}_{p}^{s+2-\frac{2}{q}-\varepsilon,\gamma+2-\frac{2}{q}-\varepsilon}(\mathbb{B})\oplus\mathbb{C}_{\omega})\\
&\cap\bigcap_{\nu\geq0,\varepsilon>0}C([\tau,T];\mathcal{H}_{p}^{\nu,\gamma+2-\frac{2}{q}-\varepsilon}(\mathbb{B})\oplus\mathbb{C}_{\omega}).
\end{split}
\end{equation}
Therefore, by \eqref{BAlg22}, \eqref{solone} and \eqref{soltwo} we conclude that
$$
\partial_{t}f=\frac{4}{n-1}u^{\frac{5-n}{n-1}}\partial_{t}u\in L^{q}([0,T];X_{0}^{s})\cap \bigcap_{\nu\geq0}L^{q}([\tau,T];X_{0}^{\nu}),
$$
and hence 
$$
f\in W^{1,q}([0,T];X_{0}^{s})\cap \bigcap_{\nu\geq0}W^{1,q}([\tau,T];X_{0}^{\nu}).
$$

Let $\mu_{j}:\Omega_{j}\subset\mathcal{B} \rightarrow\mathbb{R}^{n+1}$, $j\in\{1,...,N\}$, $N\in\mathbb{N}\backslash\{0\}$, be a covering of $\mathcal{B}$ by coordinate charts, let $\{\omega_{j}\}_{j\in\{1,...,N\}}$ be a subordinated partition of unity and let $\{y^{1},...,y^{n+1}\}$ be local coordinates in $\Omega_{j}$. On the collar part $[0,1)\times\partial\mathcal{B}$ we have $\mu_{j}:\Omega_{j}\subset\mathcal{B} \rightarrow[0,\infty)\times\mathbb{R}^{n}$ and we choose $y^{1}=x$. Denote $\{g^{ij}\}_{i,j\in\{1,...,n+1\}}=\{g_{ij}\}_{i,j\in\{1,...,n+1\}}^{-1}$ and recall the identity 
\begin{gather}\label{deltaident}
\Delta u^{m}=mu^{m-1}\Delta u+m(m-1)u^{m-2}\langle\nabla u,\nabla u\rangle_{g_{0}}, \quad m\in\mathbb{R},
\end{gather}
where, in local coordinates,
$$
\nabla u=\sum_{i,j=1}^{n+1}g^{ij}\frac{\partial u}{\partial y^{i}}\frac{\partial }{\partial y^{j}},
$$
and $\langle\cdot,\cdot\rangle_{g_{0}}$ is the Riemannian scalar product induced by $g_{0}$. In particular, on $(0,1)\times\partial\mathcal{B}$ we have
\begin{gather}\label{hoosts}
\langle\nabla u,\nabla v\rangle_{g_{0}}=\frac{1}{x^{2}}\Big((x\partial_{x}u)(x\partial_{x}v)+\sum_{i,j=1}^{n}h^{ij}(x)\frac{\partial u}{\partial y^{i+1}}\frac{\partial v}{\partial y^{j+1}}\Big),
\end{gather}
where $\{h^{ij}\}_{i,j\in\{1,...,n\}}=\{h_{ij}\}_{i,j\in\{1,...,n\}}^{-1}$.

By \eqref{YF1}-\eqref{YF2} and \eqref{deltaident} the conformal factor $f$ satisfies 
\begin{eqnarray}\nonumber
\lefteqn{u^{\frac{9-n}{n-1}}(t)u'(t)-\frac{n(n-1)}{4}\Delta f(t)}\\\label{YFf1}
&=&\frac{n(n-5)}{n-1}u^{2\frac{3-n}{n-1}}(t)\langle\nabla u(t),\nabla u(t)\rangle_{g_{0}}-\frac{n-1}{4}u^{\frac{4}{n-1}}(t)R_{g_{0}}, \quad t\in(0,T),\\\label{YFf2}
f(0)&=&1.
\end{eqnarray}
Due to \eqref{BAlg22}, \eqref{solone}, \eqref{soltwo} and \eqref{ert}, we have that
\begin{gather}\label{popsops}
u^{\frac{9-n}{n-1}}\partial_{t}u\in L^{q}([0,T];\mathcal{H}_{p}^{s,\gamma}(\mathbb{B}))\cap \bigcap_{\nu\geq0} L^{q}([\tau,T];\mathcal{H}_{p}^{\nu,\gamma}(\mathbb{B})).
\end{gather}

Let $\widetilde{x}$ be a $C^{\infty}(\mathbb{B}^{\circ})$ function such that $\widetilde{x}=x$ in local coordinates on $(0,1)\times\partial\mathcal{B}$ and $\widetilde{x}\geq\frac{1}{2}$ on $\mathcal{B}\backslash((0,1)\times\partial\mathcal{B})$. In each $\Omega_{j}$, $j\in\{1,...,N\}$, denote by $\partial_{z^{i}}$ the local derivative $\partial_{y^{i}}$ with the convention that when we are on the collar part $(0,1)\times\partial\mathcal{B}$ we denote $\partial_{z^{1}}=\partial_{x}$ and $\partial_{z^{i}}=\frac{1}{x}\partial_{y^{i}}$, $i\in\{2,...,n+1\}$. We have that 
\begin{gather}\label{edtris}
\omega_{j}^{\frac{1}{2}}\widetilde{x}^{-1}\partial_{z^{i}}u \in L^{q}([0,T];\mathcal{H}_{p}^{s+1,\gamma}(\mathbb{B}))\cap \bigcap_{\nu\geq0} L^{q}([\tau,T];\mathcal{H}_{p}^{\nu,\gamma}(\mathbb{B})),
\end{gather}
and, by \eqref{ert}, that
\begin{equation}\label{ert22}
\begin{split}
\omega_{j}^{\frac{1}{2}}\widetilde{x}\partial_{z^{i}}u \in &\bigcap_{\varepsilon>0}C([0,T];\mathcal{H}_{p}^{s+1-\frac{2}{q}-\varepsilon,\gamma+2-\frac{2}{q}-\varepsilon}(\mathbb{B}))\\
&\cap\bigcap_{\nu\geq0,\varepsilon>0}C([\tau,T];\mathcal{H}_{p}^{\nu,\gamma+2-\frac{2}{q}-\varepsilon}(\mathbb{B})).
\end{split}
\end{equation}
Therefore, by \eqref{BAlg22}, \eqref{ert}, \eqref{hoosts}, \eqref{edtris} and \eqref{ert22} we deduce that
\begin{gather}\label{erdgan}
u^{2\frac{3-n}{n-1}}\langle\nabla u,\nabla u\rangle_{g_{0}}\in L^{q}([0,T];\mathcal{H}_{p}^{s,\gamma}(\mathbb{B}))\cap \bigcap_{\nu\geq0} L^{q}([\tau,T];\mathcal{H}_{p}^{\nu,\gamma}(\mathbb{B})).
\end{gather}

Finally, by \eqref{BAlg22} and \eqref{ert} we have that
\begin{gather}\label{erdgan33}
u^{\frac{4}{n-1}}R_{g_{0}}\in L^{q}([0,T];\mathcal{H}_{p}^{s,\gamma}(\mathbb{B}))\cap \bigcap_{\nu\geq0} L^{q}([\tau,T];\mathcal{H}_{p}^{\nu,\gamma}(\mathbb{B})).
\end{gather}
Hence, from \eqref{YFf1}, \eqref{popsops}, \eqref{erdgan} and \eqref{erdgan33} we deduce that
\begin{gather}\label{kaops}
\Delta f\in L^{q}([0,T];\mathcal{H}_{p}^{s,\gamma}(\mathbb{B}))\cap \bigcap_{\nu\geq0} L^{q}([\tau,T];\mathcal{H}_{p}^{\nu,\gamma}(\mathbb{B})).
\end{gather}

Since $X_{1}^{\nu}$, $\nu\geq s$, is a Banach algebra due to \eqref{BAlg}, by \eqref{solone} and \eqref{soltwo} we have that $f(t)\in X_{1}^{\nu}$ for almost all $t\in I_{\nu}$, where $I_{\nu}=(0,T)$ when $\nu=s$ and $I_{\nu}=(\tau,T)$ when $\nu>s$. Moreover, for each $\nu\geq s$ we have that
$$
\|f(t)\|_{X_{1}^{\nu}}\leq C_{8}(\|f(t)\|_{X_{0}^{\nu}}+\|\Delta f(t)\|_{X_{0}^{\nu}}) \quad \text{for almost all} \quad t\in I_{\nu},
$$
and for certain $C_{8}>0$. Therefore, by \eqref{ert}, \eqref{kaops} and Minkowski inequality we conclude that
\begin{gather*}
 f\in L^{q}([0,T];\mathcal{H}_{p}^{s+2,\gamma+2}(\mathbb{B})\oplus\mathbb{C}_{\omega})\cap \bigcap_{\nu\geq0} L^{q}([\tau,T];\mathcal{H}_{p}^{\nu,\gamma+2}(\mathbb{B})\oplus \mathbb{C}_{\omega}).
\end{gather*}
The continuity follows from \eqref{interpemb} and \eqref{embd11}. \mbox{\ } \hfill $\square$

\subsection*{Proof of Corollary 1.2}
The result follows by Theorem \ref{tshortyf} and \cite[Corollary 2.9]{RS2}. \mbox{\ } \hfill $\square$

\end{document}